\documentclass{amsart}

\usepackage{fancyhdr}
\usepackage{lastpage}
\usepackage{yhmath}

\newcommand{\bdis}{\begin{displaymath}}
\newcommand{\edis}{\end{displaymath}}
\newcommand{\be}{\begin{equation}}
\newcommand{\ee}{\end{equation}}

\newcommand{\mbb}{\mathbb}
\newcommand{\mcal}{\mathcal}

\newcommand{\vp}{\varphi}
\newcommand{\vt}{\vartheta}

\newcommand{\zf}{\zeta\left(\frac{1}{2}+it\right)}

\DeclareMathOperator*{\ssum}{\sum\sum}

\theoremstyle{definition}

\theoremstyle{remark}
\newtheorem{remark}[]{Remark}

\newtheorem*{mydef1}{{\bf Theorem}}

\newtheorem*{mydef11}{{\bf Theorem 1}}

\newtheorem*{mydef12}{{\bf Theorem 2}}

\newtheorem*{mydef41}{{\bf Corollary 1}}

\newtheorem*{mydef42}{{\bf Corollary 2}}
\newtheorem*{mydef43}{{\bf Corollary 3}}
\newtheorem*{mydef44}{{\bf Corollary 4}}

\newtheorem*{mydef51}{{\bf Lemma 1}}

\newtheorem*{mydef52}{{\bf Lemma 2}}

\newtheorem*{mydef53}{{\bf Lemma 3}}

\newtheorem*{mydef54}{{\bf Lemma 4}}

\newtheorem*{mydef55}{{\bf Lemma 5}}

\newtheorem*{mydef56}{{\bf Lemma 6}}

\newtheorem*{mydef57}{{\bf Lemma 7}}

\numberwithin{equation}{section}

\begin{document}

\title{Lindel\" of hypothesis and the order of the mean-value of $|\zeta(s)|^{2k-1}$ in the critical strip}

\author{Jan Moser}

\address{Department of Mathematical Analysis and Numerical Mathematics, Comenius University, Mlynska Dolina M105, 842 48 Bratislava, SLOVAKIA}

\email{jan.mozer@fmph.uniba.sk}


\begin{abstract}
The main subject of this paper is the mean-value of the function $|\zeta(s)|^{2k-1}$ in the critical strip.
On Lindel\" of hypothesis we give a solution to this question for some class of disconnected sets. This paper is English
version of our paper \cite{5}.
\end{abstract}

\maketitle

\section{Introduction}

\subsection{}

E.C. Titchmarsh had began with the study of the mean-value of the function
\bdis
\left|\zeta\left(\sigma+it\right)\right|^\omega,\ \frac 12<\sigma\leq 1,\ 0<\omega ,
\edis
where $\omega$ is non-integer number, \cite{6} (comp. \cite{2}, p. 278). Next, Ingham and Davenport have obtained
the following result (see \cite{1}, \cite{2}, comp. \cite{7}, pp. 132, 133)
\be \label{1.1}
\frac 1T \int_1^T \left|\zeta\left(\sigma+it\right)\right|^{2\omega}{\rm d}t=
\sum_{n=1}^\infty \frac{d^2_{\omega}(n)}{n^{2\sigma}}+\mcal{O}(1),\ \omega\in (0,2],\ T\to\infty.
\ee
Let us remind that:
\begin{itemize}
 \item[(a)] for $\omega\in\mbb{N}$ the symbol $d_\omega(n)$ denotes the number of decompositions of $n$ into
 $\omega$-factors ,
 \item[(b)] in the case $\omega$ is not an integer, we define $d_\omega(n)$ as the coefficient of $n^{-s}$ in the
 Dirichlet series for the function $\zeta^\omega(s)$ converging for all $\sigma>1$.
\end{itemize}

\subsection{}

Next, for
\bdis
\omega=\frac 12,\frac 32
\edis
it follows from (\ref{1.1})  that the orders of mean-values
\bdis
\frac 1T \int_1^T \left|\zeta\left(\sigma+it\right)\right|{\rm d}t,\
\frac 1T \int_1^T \left|\zeta\left(\sigma+it\right)\right|^3{\rm d}t
\edis
are determined. But a question about the order of mean-value of
\bdis
|\zeta(\sigma+it)|^{2l+1},\ l=2,3,\dots,\ \frac 12<\sigma<1
\edis
remains open.

In this paper we give a solution to this open question on the assumption of truth of the Lindel\" of
hypothesis for some infinite class of disconnected sets. In a particular case we obtain the following result: on
Lindel\" of hypothesis we have
\be \label{1.2}
\begin{split}
 & 1-|o(1)|<\frac 1H\int_T^{T+H}|\zeta(\sigma+it)|^{2k-1}{\rm d}t < \\
 & < \sqrt{F(\sigma,2k-1)}+|o(1)|,\quad H=T^\epsilon,\ k=1,2,\dots,\ 0<\epsilon ,
\end{split}
\ee
where
\be \label{1.3}
F(\sigma,\omega)=\sum_{n=1}^\infty \frac{d^2_{\omega}(n)}{n^{2\sigma}},
\ee
and $\epsilon$ is an arbitrarily small number.

The proof of our main result is based on our method (see \cite{4}) for the proof of new mean-value theorem for the
Riemann zeta-function
\bdis
Z(t)=e^{i\vt(t)}\zf
\edis
with respect of two infinite classes of disconnected sets.

\section{Main formulas}

We use the following formula: on Lindel\" of hypothesis
\be \label{2.1}
\begin{split}
 & \zeta^k(s)=\sum_{n\leq t^\delta}\frac{d_k(n)}{n^s}+\mcal{O}(t^{-\lambda}),\ \lambda=\lambda(k,\delta,\sigma)>0, \\
 & s=\sigma+it,\ \frac 12<\sigma<1,\ t>0
\end{split}
\ee
(see \cite{7}, p. 277) for every natural number $k$, where $\delta$ is any given positive number less than $1$. Let us
remind that
\be \label{2.2}
d_k(n)=\mcal{O}(n^\eta) ,
\ee
where $0<\eta$ is an arbitrarily small number. Of course, (see (\ref{2.1}), (\ref{2.2}))
\be \label{2.3}
\begin{split}
 & \zeta^k(s)=\mcal{O}\left(\sum_{n\leq t^\delta}d_k(n)n^{-\sigma}\right)=
 \mcal{O}\left( t^{\delta\eta+\delta(1-\sigma)}\right)= \\
 & = \mcal{O}\left( t^{(n+1/2)\delta}\right).
\end{split}
\ee
Let
\be \label{2.4}
t\in [T,T+H],\ H=T^{\epsilon},\ 2\delta n+2\delta<\epsilon.
\ee
Since
\bdis
\sum_{T^\delta\leq n\leq (T+H)^\delta}1=\mcal{O}(T^{\delta+\epsilon-1})
\edis
then
\be \label{2.5}
\begin{split}
 & \sum_{T^\delta\leq n\leq t^{\delta}}\frac{d_k(n)}{n^s}=
 \mcal{O}\left( T^{\delta\eta+\delta\sigma}\cdot \sum_{T^\delta\leq n\leq (T+H)^\delta} 1 \right)= \\
 & = \mcal{O}(T^{\delta\eta-\delta\sigma+\delta+\epsilon-1})=\mcal{O}(T^{-\lambda_1}),
\end{split}
\ee
where
\be \label{2.6}
\lambda_1=1-\delta-\epsilon+\delta\sigma-\delta\eta>0 ,
\ee
(of course, for sufficiently small $\epsilon$ the inequality (\ref{2.6}) holds true). Next, for
\be \label{2.7}
\lambda_2=\lambda_2(k,\delta,\sigma,\epsilon,\eta)=\min \{\lambda,\lambda_1\}>0
\ee
the following formula (see (\ref{2.1}), (\ref{2.5}) -- (\ref{2.7}))
\be \label{2.8}
\zeta^k(s)=\sum_{n<T^\delta}\frac{d_k(n)}{n^s}+\mcal{O}(T^{-\lambda_2}),\ t\in [T,T+H]
\ee
holds true. Since
\be \label{2.9}
\zeta^k(s)=U_k(\sigma,t)+iV_k(\sigma,t)
\ee
then - on Lindel\" of hypothesis - we obtain from (\ref{2.8}) the following main formula
\be \label{2.10}
\begin{split}
 & U_k(\sigma,t)=1+\sum_{2\leq n < T^\delta}\frac{d_k(n)}{n^\sigma}\cos(t\ln n)+\mcal{O}(T^{-\lambda_2}), \\
 & V_k(\sigma,t)=-\sum_{2\leq n < T^\delta}\frac{d_k(n)}{n^\sigma}\sin(t\ln n)+\mcal{O}(T^{-\lambda_2}), \\
 & t\in [T,T+H].
\end{split}
\ee

\section{The first class of lemmas}

Let us denote by
\bdis
\{ t_\nu(\tau)\}
\edis
an infinite set of sequences that we defined (see \cite{4}, (1)) by the condition
\be \label{3.1}
\vt[t_\nu(\tau)]=\pi\nu+\tau,\ \nu=1,\dots,\ \tau\in [-\pi,\pi],
\ee
of course,
\bdis
t_\nu(0)=t_\nu,
\edis
where (see \cite{7}, pp. 220, 329)
\be \label{3.2}
\begin{split}
 & \vt(t)=-\frac t2\ln\pi+\mbox{Im}\ln\Gamma\left(\frac 12+i\frac t2\right), \\
 & \vt'(t)=\frac 12\ln\frac{t}{2\pi}+\mcal{O}\left(\frac 1t\right), \\
 & \vt''(t)\sim \frac{1}{2t}.
\end{split}
\ee

\subsection{}

The following lemma holds true.

\begin{mydef51}
If
\bdis
2\leq m,n < T^{\delta}
\edis
then
\be \label{3.3}
\sum_{T\leq t_\nu\leq T+H}\cos\{ t_\nu(\tau)\ln n\}=\mcal{O}\left(\frac{\ln T}{\ln n}\right),
\ee
\be \label{3.4}
\sum_{T\leq t_\nu\leq T+H}\cos\{ t_\nu(\tau)\ln (mn)\}=\mcal{O}\left(\frac{\ln T}{\ln (mn)}\right),
\ee
\be \label{3.5}
\sum_{T\leq t_\nu\leq T+H}\cos\left\{ t_\nu(\tau)\ln \frac mn\right\}=
\mcal{O}\left(\frac{\ln T}{\ln \frac mn}\right), \ m>n,
\ee
where the $\mcal{O}$-estimates are valid uniformly for $\tau\in [-\pi,\pi]$.
\end{mydef51}

\begin{proof}
We use the van der Corput's method. Let (see (\ref{3.3}))
\bdis
\vp_1(\nu)=\frac{1}{2\pi}t_\nu(\tau)\ln n.
\edis
Next, (see (\ref{2.4}), (\ref{3.1}), (\ref{3.2}))
\bdis
\begin{split}
 & \vp_1'(\nu)=\frac{\ln n}{2\vt'[t_\nu(\tau)]}, \\
 & \vp_1''(\nu)=-\frac{\pi\ln n}{2\{ \vt'[t_\nu(\tau)]\}^2}\vt''\{ t_\nu(\tau)\}<0, \\
 & 0< A\frac{\ln n}{\ln T}\leq \vp_1'(\nu)=
 \frac{\ln n}{\ln\frac{t_\nu(\tau)}{2\pi}+\mcal{O}(\frac 1t)}=\frac{\ln n}{\ln\frac{T}{2\pi}+\mcal{O}(\frac HT)}< \\
 & < \delta\frac{\ln T}{\ln\frac{T}{2\pi}+\mcal{O}(\frac HT)}<\frac 14,
\end{split}
\edis
($A>0$, sice $\delta$ may be sufficiently small). Hence, (see \cite{7}, p. 65 and p. 61, Lemma 4.2)
\bdis
\begin{split}
 & \sum_{T\leq t_\nu\leq T+H}\cos\{ t_\nu(\tau)\ln n\}= \\
 & = \int_{T\leq t_x\leq T+H} \cos\{ 2\pi\vp_1(x)\}{\rm d}x+\mcal{O}(1)=
 \mcal{O}\left(\frac{\ln T}{\ln n}\right),
\end{split}
\edis
i.e. the estimate (\ref{3.3}) holds true. The estimates (\ref{3.4}) and (\ref{3.5}) follow by the similar way.
\end{proof}

\subsection{}

The following lemma holds true.

\begin{mydef52}
On Lindel\" of hypothesis we have
\be \label{3.6}
\sum_{T\leq t_\nu\leq T+H} U_k[\sigma,t_\nu(\tau)]=\frac{1}{2\pi}H\ln\frac{T}{2\pi}+\mcal{O}(H).
\ee
\end{mydef52}

\begin{proof}
Let us remind that
\be \label{3.7}
\sum_{T\leq t_\nu\leq T+H} 1=\frac{1}{2\pi}H\ln\frac{T}{2\pi}+\mcal{O}(1),
\ee
(see \cite{3}, (23)). Next, (see (\ref{2.10}), (\ref{3.7}))
\be \label{3.8}
\begin{split}
 & \sum_{T\leq t_\nu\leq T+H} U_k[\sigma,t_\nu(\tau)]=\frac{1}{2\pi}H\ln\frac{T}{2\pi}+\mcal{O}(1)+\\
 & + \mcal{O}(T^{-\lambda_2}H\ln T)+\sum_{2\leq n<T^\delta}\frac{d_k(n)}{n^\sigma}\cdot
 \sum_{T\leq t_nu\leq T+H}\cos\{ t_\nu(\tau)\ln n\}= \\
 & = \frac{1}{2\pi}H\ln T+\mcal{O}(1)+\mcal{O}(T^{-\lambda_2}H\ln T)+w_1.
\end{split}
\ee
Since (see (\ref{2.2}), (\ref{2.4}), (\ref{3.3}))
\bdis
\begin{split}
 & w_1=\mcal{O}\left(T^{\delta\eta}\ln T\sum_{2\leq n\leq T^\delta}\frac{1}{\sqrt{n}\ln n}\right)= \\
 & = \mcal{O}\left\{ T^{\delta\eta}\ln T\left( \sum_{2\leq n<T^{\delta/2}} \ + \
 \sum_{T^{\delta/2}\leq n<T^\delta}\right)\frac{1}{\sqrt{n}\ln n}\right\}= \\
 & = \mcal{O}(T^{\delta\eta+\delta/2})=\mcal{O}(H),
\end{split}
\edis
then from (\ref{3.8}) the formula (\ref{3.6}) follows.
\end{proof}

\section{Theorem 1}

\subsection{}

Next, we define the following class of disconnected sets (comp. \cite{4}, (3)):
\be \label{4.1}
G(x)=\bigcup_{T\leq t_\nu\leq T+H}\{ t:\ t_\nu(-x)<t<t_\nu(x)\},\ 0<x\leq \frac{\pi}{2}.
\ee
Let us remind that (see \cite{4}, (7))
\be \label{4.2}
\begin{split}
 & t_\nu(x)-t_\nu(-x)=\frac{4x}{\ln\frac{T}{2\pi}}+\mcal{O}\left(\frac{xH}{T\ln^2T}\right), \\
 & t_\nu(-x),t_\nu(x)\in [T,T+H].
\end{split}
\ee
Of course,
\be \label{4.3}
m\{ G(x)\}=\frac{2x}{\pi}H+\mcal{O}(x),
\ee
(see (\ref{3.7}), (\ref{4.2})), where $m\{ G(x)\}$ stands for the measure of $G(x)$.

\subsection{}

The following theorem holds true.

\begin{mydef11}
On Lindel\" of hypothesis
\be \label{4.4}
\int_{G(x)} U_k(\sigma,t){\rm d}t=\frac{2x}{\pi}H+o\left(\frac{xH}{\ln T}\right).
\ee
\end{mydef11}

First of all, we obtain from (\ref{4.4}) by (\ref{4.3}) the following

\begin{mydef41}
On Lindel\" of hypothesis
\be \label{4.5}
\frac{1}{m\{ G(x)\}}\int_{G(x)} U_k(\sigma,t){\rm d}t=1+o\left(\frac{1}{\ln T}\right).
\ee
\end{mydef41}

Next, we obtain from (\ref{4.4}) the following

\begin{mydef42}
On Lindel\" of hypothesis
\be \label{4.6}
\int_{G(x)}|U_k(\sigma,t)|{\rm d}t\geq \frac{2xH}{\pi}\{ 1-|o(1)|\}.
\ee
\end{mydef42}

Since (see (\ref{2.9}))
\bdis
|\zeta(s)|^{2k-1}=\sqrt{U^2_{2k-1}+V^2_{2k-1}}\geq |U_{2k-1}|
\edis
then we obtain from(\ref{4.6}), $k\longrightarrow 2k-1$, the following

\begin{mydef43}
On Lindel\" of hypothesis
\be \label{4.7}
\int_{G(x)}|\zeta(\sigma+it)|^{2k-1}{\rm d}t\geq \frac{2xH}{\pi}\{ 1-|o(1)|\}.
\ee
\end{mydef43}

\subsection{}

In this part we shall give the

\begin{proof}
of the Theorem 1. Since (see (\ref{3.1}), (\ref{3.2}))
\bdis
\begin{split}
& \left(\frac{{\rm d}t_\nu(\tau)}{{\rm d}\tau}\right)^{-1}=\ln P_0+\mcal{O}\left(\frac HT\right),\\
& t_\nu(\tau)\in [T,T+H],\ P_0=\sqrt{\frac{T}{2\pi}},
\end{split}
\edis
then we obtain by using of the substitution
\bdis
t=t_\nu(\tau) ,
\edis
and by estimates (\ref{2.3}) that
\bdis
\begin{split}
 & \int_{-\pi}^\pi U_k[\sigma,t_\nu(\tau)]{\rm d}\tau=\int_{-\pi}^\pi U_k[\sigma,t_\nu(\tau)]
 \left(\frac{{\rm d}t_\nu(\tau)}{{\rm d}\tau}\right)^{-1}\cdot \frac{{\rm d}t_\nu(\tau)}{{\rm d}\tau}
 {\rm d}\tau= \\
 & = \ln P_0\int_{-\pi}^\pi U_k[\sigma,t_\nu(\tau)]\frac{{\rm d}t_\nu(\tau)}{{\rm d}\tau}{\rm d}\tau+ \\
 & + \mcal{O}\left( x\max\{|\zeta^k|\}\frac HT\max\left\{\frac{{\rm d}t_\nu(\tau)}{{\rm d}\tau}\right\}\right)= \\
 & = \ln P_0\int_{t_\nu(-x)}^{t_\nu(x)}U_k(\sigma,t){\rm d}t+
 \mcal{O}\left( x\frac{T^{\delta\eta+\delta/2+\epsilon-1}}{\ln T}\right),
\end{split}
\edis
where the $\max$ is taken with respect to the segment $[T,T+H]$. Consequently, (see (\ref{2.3}),
(\ref{3.7}), (\ref{4.1}) and (\ref{4.2}))
\bdis
\begin{split}
 & \sum_{T\leq t_\nu\leq T+H}\int_{-\pi}^\pi U_k[\sigma,t_\nu(\tau)]{\rm d}\tau= \\
 & = \ln P_0 \int_{G(x)}U_k(\sigma,t){\rm d}t+\mcal{O}(xT^{\delta\eta+\delta/2+2\epsilon-1})+ \\
 & + \mcal{O}\left(\frac{xT^{(\eta+1/2)\delta}}{\ln T}\right).
\end{split}
\edis
Now, the integration (\ref{3.6}) by
\bdis
\tau\in [-\pi,\pi]
\edis
gives the formula
\bdis
\begin{split}
 & \ln P_0\int_{G(x)} U_k(\sigma,t){\rm d}t+\mcal{O}(xT^{\delta\eta+\delta/2+2\epsilon-1})= \\
 & = \frac x\pi H\ln\frac{T}{2\pi}+\mcal{O}(xT^{\delta\eta+\delta/2})
\end{split}
\edis
and from this by (\ref{2.4}) the formula (\ref{4.4}) follows immediately (here $\epsilon$ is
arbitrarily small number).
\end{proof}

\section{The second class of lemmas}

\subsection{}

Let
\be \label{5.1}
\begin{split}
 & S_1(t)=\sum_{2\leq n\leq T^\delta}\frac{d_k(n)}{n^\sigma}\cos(t\ln n), \\
 & w_2(t)=\{ S_1(t)\}^2.
\end{split}
\ee
The following lemma holds true.

\begin{mydef53}
\be \label{5.2}
\begin{split}
 & \sum_{T\leq t_\nu\leq T+H} w_2[t_\nu(\tau)]= \\
 & = \{ F(\sigma,k)-1\}\cdot \frac{1}{4\pi} H\ln\frac{T}{2\pi}+o(H),
\end{split}
\ee
(on $F(\sigma,k)$ see (\ref{1.3})).
\end{mydef53}

\begin{proof}
First of all we have
\be \label{5.3}
\begin{split}
 & w_2(t)=\sum_m\sum_n \frac{d_k(m)d_k(n)}{(mn)^\sigma}\cos(t\ln m)\cos(t\ln n)= \\
 & = \frac 12\sum_m\sum_n\frac{d_k(m)d_k(n)}{(mn)^\sigma}\cos\{t\ln(mn)\}+\\
 & + \ssum_{n<m}\frac{d_k(m)d_k(n)}{(mn)^\sigma}\cos\left(t\ln\frac mn\right)+
 \frac 12\sum_n\frac{d_k^2(n)}{n^{2\sigma}}= \\
 & = w_{21}(t)+w_{22}(t)+w_{23}(t).
\end{split}
\ee
Now we have:

by (\ref{2.2}), (\ref{2.4}) and (\ref{3.4})
\be \label{5.4}
\begin{split}
 & \sum_{T\leq t_\nu\leq T+H}w_{21}[t_\nu(\tau)]=
 \mcal{O}\left( T^{\delta\eta}\ln T\cdot \ssum_{2\leq m,n<T^\delta}\frac{1}{\sqrt{mn}\ln(mn)}\right)= \\
 & = \mcal{O}(T^{2\delta\eta+\delta}\ln T)=o(H);
\end{split}
\ee
by (\ref{2.2}), (\ref{2.4}) and (\ref{3.5}) and by \cite{7}, p. 116, Lemma, ($T\longrightarrow T^\delta$),
\be \label{5.5}
\begin{split}
 & \sum_{T\leq t_\nu\leq T+H}w_{22}[t_\nu(\tau)]=
 \mcal{O}\left( T^{\delta\eta}\ln T\cdot \ssum_{2\leq n<m<T^\delta}\frac{1}{\sqrt{mn}\ln\frac mn}\right)= \\
 & = \mcal{O}(T^{2\delta\eta+\delta}\ln^2T)=o(H);
\end{split}
\ee
by (\ref{1.3}), $\omega\longrightarrow k$, and by (\ref{2.2})
\be\label{5.6}
\begin{split}
 & 2_{23}=\frac 12\sum_{n=2}^\infty \frac{d_k^2(n)}{n^{2\sigma}}-\frac 12\sum_{n\geq T^{\delta}}
 \frac{d_k^2(n)}{n^{2\sigma}}= \\
 & = \frac 12\{ F(\sigma,k)-1\}+\mcal{O}\left(\int_{T^\delta}^\infty x^{\eta-2\sigma}{\rm d}x\right)= \\
 & = \frac 12\{ F(\sigma,k)-1\}+\mcal{O}(T^{-\delta(2\sigma-1-\eta)});
\end{split}
\ee
(of course, $2\sigma-1-\eta>0$ since $\eta$ is arbitrarily small). Finally, by (\ref{2.4}), (\ref{3.7})
and (\ref{5.6}) we obtain
\be \label{5.7}
\sum_{T\leq t_\nu\leq T+H}w_{23}=\{ F(\sigma,k)-1\}\frac{1}{4\pi}H\ln\frac{T}{2\pi}+o(H).
\ee
Hence, from (\ref{5.3}) by (\ref{5.4}) -- (\ref{5.7}) the formula (\ref{5.2}) follows.
\end{proof}

Next, the following lemma holds true.

\begin{mydef54}
On Lindel\" of hypothesis
\be \label{5.8}
\sum_{T\leq t_\nu\leq T+H}U_k^2[\sigma,t_\nu(\tau)]=\{ F(\sigma,k)-1\}\frac{1}{2\pi}G\ln\frac{T}{2\pi}+
o(H).
\ee
\end{mydef54}

\begin{proof}
Since (see (\ref{2.10}), (\ref{5.1}))
\bdis
U_k(\sigma,t)=1+S_1+\mcal{O}(T^{-\lambda_2}) ,
\edis
then
\be \label{5.9}
U_k^2(\sigma,t)=1+w_2+2S_1+\mcal{O}(|S_1|T^{-\lambda_2})+\mcal{O}(T^{-2\lambda_2}).
\ee
Now we have:

by (\ref{2.4})
\bdis
\sum_{T\leq t_\nu\leq T+H}S_1[t_\nu(\tau)]=\mcal{O}(T^{\delta\eta+\delta/2})=o(H);
\edis
by(\ref{2.4}), (\ref{3.7}) and (\ref{5.2})
\bdis
\begin{split}
 & \sum_{T\leq t_\nu\leq T+H}|S_1|T^{-\lambda_2}= \\
 & = \mcal{O}\left\{ T^{-\lambda_2}\sqrt{H\ln T}
 \left(\sum_{T\leq t_\nu\leq T+H}w_2[t_\nu(\tau)]\right)^{1/2}\right\}=\\
 & = \mcal{O}(T^{-\lambda_2}H\ln T)=o(H).
\end{split}
\edis
Consequently, from (\ref{5.9}) by (\ref{3.7}) the formula (\ref{5.8}) follows.
\end{proof}

\subsection{}

Let
\be \label{5.10} \begin{split}
& S_2(t)=\sum_{2\leq n<T^\delta}\frac{d_k(n)}{n^\sigma}\sin(t\ln n), \\
& w_3(t)=\{ S_2(t)\}^2. \end{split}
\ee
The following lemma holds true.

\begin{mydef55}
\be \label{5.11}
\sum_{T\leq t_\nu\leq T+H} w_3[t_\nu(\tau)]=\{ F(\sigma,k)-1\}\frac{1}{4\pi}H\ln\frac{T}{2\pi}+o(H).
\ee
\end{mydef55}

\begin{proof}
Since (comp. (\ref{5.3}))
\be \label{5.12}
\begin{split}
 & w_3(t)=\ssum_{m,n}\frac{d_k(m)d_k(n)}{(mn)^\sigma}\sin(t\ln m)\sin(t\ln n)= \\
 & = -\frac 12\ssum_{m,n}\frac{d_k(m)d_k(n)}{(mn)^\sigma}\cos\{ t\ln(mn)\}+ \\
 & + \ssum_{n<m} \frac{d_k(m)d_k(n)}{(mn)^\sigma}\cos\left( t\ln\frac mn\right)+ \\
 & + \frac 12 \sum_n \frac{d_k^2(n)}{n^{2\sigma}}=w_{31}(t)+w_{32}(t)+w_{33}(t),
\end{split}
\ee
then we obtain by the way (\ref{5.3}) -- (\ref{5.7}) our formula (\ref{5.11}).
\end{proof}

Next, the following lemma holds true

\begin{mydef56}
On Lindel\" of hypothesis
\be \label{5.13}
\begin{split}
& \sum_{T\leq t_\nu\leq T+H} V_k^2[\sigma,t_\nu(\tau)]= \\
& = \{ F(\sigma,k)-1\}\frac{1}{4\pi}H\ln\frac{T}{2\pi}+o(H).
\end{split}
\ee
\end{mydef56}

\begin{proof}
Since (see (\ref{2.10}), (\ref{5.10}))
\bdis
V_k(\sigma,t)=-S_2+\mcal{O}(T^{-\lambda_2}),
\edis
then
\be \label{5.14}
V_k^2(\sigma,t)=w_3+\mcal{O}(T^{-\lambda_2}|S_2|)+\mcal{O}(T^{-2\lambda_2}).
\ee
Consequently, the proof may be finished in the same way as it was done in the case of our Lemma 4 (comp.
(\ref{5.12}), (\ref{5.14})).
\end{proof}

Since (see (\ref{2.9}))
\bdis
|\zeta(s)|^{2k}=U_k^2+V_k^2 ,
\edis
then by (\ref{5.8}), (\ref{5.11}) we obtain the following.

\begin{mydef57}
On Lindel\" of hypothesis
\be \label{5.14}
\sum_{T\leq t_\nu\leq T+H} |\zeta[\sigma_it_\nu(\tau)]|^{2k}=\frac{1}{2\pi}F(\sigma,k)H\ln\frac{T}{2\pi}+o(H).
\ee
\end{mydef57}

\section{Theorem 2 and main Theorem}

Now we obtain from (\ref{5.14}) by the way very similar to than one used in the proof of the Theorem 1, the following.

\begin{mydef12}
On Lindel\" of hypothesis
\be \label{6.1}
\int_{G(x)}|\zeta(\sigma+it)|^{2k}{\rm d}t=\frac{2x}{\pi}F(\sigma,t)H+o\left(\frac{xH}{\ln T}\right).
\ee
\end{mydef12}

Further, from (\ref{6.1}) we obtain

\begin{mydef44}
On Lindel\" of hypothesis
\be \label{6.2}
\int_{G(x)}|\zeta(\sigma+it)|^{2k-1}{\rm d}t<\frac{2xH}{\pi}\sqrt{F(\sigma,2k-1)}\cdot \{1+|o(1)|\}.
\ee
\end{mydef44}

Indeed, by (\ref{4.3}), (\ref{6.1}) we have
\bdis
\begin{split}
& \int_{G(x)}|\zeta(\sigma+it)|^{2k}{\rm d}t< \\
& < \sqrt{m\{ G(x)\}}\left( \int_{G(x)}|\zeta(\sigma+it)|^{4k-2}{\rm d}t\right)^{1/2}< \\
& < \frac{2xH}{\pi}\sqrt{F(\sigma,2k-1)}\cdot\{ 1+|o(1)|\}.
\end{split}
\edis

Finally, from (\ref{4.7}), (\ref{6.2}) we obtain our main result:

\begin{mydef1}
On Lindel\" of hypothesis
\be \label{6.3}
\begin{split}
& 1-|o(1)|<\frac{1}{m\{ G(x)\}}\int_{G(x)}|\zeta(\sigma+it)|^{2k-1}{\rm d}t< \\
& < \sqrt{F(\sigma,2k-1)}+|o(1)|.
\end{split}
\ee
\end{mydef1}

\begin{remark}
The question about the order of the mean-value of the function
\bdis
|\zeta(\sigma+it)|^{2k-1},\ k=1,2,\dots
\edis
defined on infinite class of disconnected sets $\{ G(x)\}$ is answered by the inequalities (\ref{6.3}).
\end{remark}

\begin{remark}
Inequalities (\ref{1.2}) follows from (\ref{6.3}) as a special case for $x=\pi/2$, (see (\ref{2.3}), (\ref{2.4}), (\ref{4.3})).
\end{remark}

\thanks{I would like to thank Michal Demetrian for helping me with the electronic version of this work.}

\end{document}